\documentstyle{article}
\begin{document}
\newtheorem{theorem}{Theorem}[section]
\newtheorem{definition}{Definition}[section]
\newtheorem{corollary}[theorem]{Corollary}
\newtheorem{lemma}[theorem]{Lemma}
\newtheorem{proposition}[theorem]{Proposition}
\newtheorem{step}[theorem]{Step}
\newtheorem{example}[theorem]{Example}
\newtheorem{remark}[theorem]{Remark}

\font\sixbb=msbm6
\font\eightbb=msbm8
\font\twelvebb=msbm10 scaled 1095
\newfam\bbfam
\textfont\bbfam=\twelvebb \scriptfont\bbfam=\eightbb
                           \scriptscriptfont\bbfam=\sixbb

\def\bb{\fam\bbfam\twelvebb}
\newcommand{\Rea}{{\bb R}}
\newcommand{\Com}{{\bb C}}
\newcommand{\Int}{{\bb Z}}
\newcommand{\Nat}{{\bb N}}
\def\calC{{\cal C}}
\def\calP{{\cal P}}
\def\calQ{{\cal Q}}
\def\calD{{\cal D}}
\newcommand{\D}{{\bb D}}
\newcommand{\E}{{\bb E}}
\newcommand{\F}{{\bb F}}
\newcommand{\PP}{{\bb P}}
\newcommand{\card}{\rm card}
\newcommand{\FF}{{\bb F}}
\newcommand{\rk}{{\rm rank}}
\newcommand{\Ker}{{\rm Ker}}
\renewcommand{\Im}{{\rm Im}}
\newcommand{\enp}{\begin{flushright} $\Box$ \end{flushright}}
\def\cD{{\mathcal{D}}}
\newcommand{\Span}{{\rm Span}}

\def\Sn{{\cal S}_n(F)}
\def\PSn{P{\cal S}_n(F)}
\def\cM{l}
\def\vp{\varphi}
\def\GL{\mbox{\rm GL}}
\def\LRa{\Leftrightarrow}
\def\bp{\noindent{\bf Proof. } }
\def\ep{\hfill $\Box$}

\def\N{{\bb N}}
\def\H{{\cal H}_n(D)}
\def\inv{\overline{\phantom{e}}}

\title{Continuous space-time transformations
\thanks{The second author was
supported by a grant from ARRS, Slovenia.}}
\author{Cl\' ement de Seguins Pazzis\footnote{Laboratoire de Math\' ematiques de Versailles,
Universit\' e de Versailles Saint-Quentin-en-Yvelines,
45 avenue des Etats-Unis, 78035 Versailles cedex, France. dsp.prof@gmail.com}, \quad
Peter \v Semrl\footnote{Faculty of Mathematics and Physics, University of Ljubljana,
        Jadranska 19, SI-1000 Ljubljana, Slovenia. peter.semrl@fmf.uni-lj.si}
        }

\date{}
\maketitle

\begin{abstract}
We prove that every continuous map acting on the four-dimensional Minkowski space and preserving light cones in one direction only
is either a Poincar\' e similarity, that is, a product of a Lorentz transformation and a dilation, or it is of a very special degenerate form.
In the presence of the continuity assumption the main tool in the proof is a basic result from the homotopy theory of spheres.
\end{abstract}
\maketitle



\section{Introduction and statement of the main result}

\subsection{The problem}

Throughout the paper $\Rea^4 = \{ (x,y,z,t)\, : \, x,y,z,t \in \Rea \}$ will be equipped with the standard topology.
At the foundations of relativity theory lies the Lorentz-Minkowski indefinite inner product defined by
$$
\left\langle (x,y,z,t) , (u,v,w,s) \right\rangle = -xu - yv -zw + ts.
$$
A fundamental result due to Alexandrov \cite{Al1,Al2,Al3} (see also \cite{Ben,Zee}) states that every bijective map $\phi : \Rea^4 \to \Rea^4$
with the property that for every pair of space-time events $(r_1, r_2) \in (\Rea^4)^2$ we have
\begin{equation}\label{ggg}
\langle r_1 - r_2 , r_1 - r_2 \rangle = 0 \iff \langle \phi (r_1 ) - \phi ( r_2 ) , \phi ( r_1 ) - \phi ( r_2 ) \rangle = 0
\end{equation}
is of the form
\begin{equation}\label{fiesa}
\phi : r \in \Rea^4 \mapsto  k\,Qr + a,
\end{equation}
where $k$ is a positive real number, $a$ is a vector in $\Rea^4$, and $Q$ is a $4\times 4$ Lorentz matrix, that is, a matrix satisfying $QMQ^t = M$, where
$$
M = \left[ \begin{matrix} { -1 & 0 & 0 & 0 \cr
0 & -1 & 0 & 0 \cr
0 & 0 & -1 & 0 \cr
0 & 0 & 0 & 1 \cr} \end{matrix} \right].
$$
For a physical interpretation of the above theorem we refer to \cite[p.691]{Pfe}. Let us just mention that for a given $r\in \Rea^4$ the set $\{ s \in \Rea^4 \, : \, \langle s - r , s - r \rangle = 0 \}$
is called the light cone with vertex $r$.

Because of its importance in the mathematical foundations of relativity theory, Alexandrov's theorem has been improved in several directions.
Instead of the four-dimensional Minkowski space one can consider Minkowski spaces of any dimension $\ge 3$ or even infinite-dimensional Minkowski spaces,
maps can be defined on domains of Minkowski spaces (then, besides Lorentz maps composed with dilations, inversions and singular double inversions satisfy property (\ref{ggg}) ),
such maps were considered with the bijectivity assumption relaxed to injectivity only or to surjectivity only, the relation (\ref{ggg}) can be replaced by a similar one
with ``$=$" replaced by either ``$\ge$", or ``$>$", etc. An interested reader can find this kind of results in \cite{BoH, ChP, Les, LeM, Pfe, PoR}  and the references therein.

When such maps are defined on the whole Minkowski space (not just on some open domain), these results can be considered as various
possible derivations of Lorentz transformations under weak mathematical properties and the proofs are often based on certain improvements of the fundamental theorem of affine geometry under weak assumptions.

In this paper we will describe the general form of continuous maps  $\phi : \Rea^4 \to \Rea^4$
with the property that for every pair of space-time events $r_1, r_2 \in \Rea^4$ we have
\begin{equation}\label{clepet}
\langle r_1 - r_2 , r_1 - r_2 \rangle = 0 \Rightarrow \langle \phi (r_1 ) - \phi ( r_2 ) , \phi ( r_1 ) - \phi ( r_2 ) \rangle = 0.
\end{equation}
We assume neither injectivity nor surjectivity and we assume that the condition $\langle r_1 - r_2 , r_1 - r_2 \rangle = 0$ is preserved in one direction only.
Still, under these very weak assumptions (besides the weak preserving property our only assumption is continuity) we get almost the same conclusion as
in Alexandrov's classical result. Namely, such maps are either Poincar\' e similarities, that is,
Lorentz transformations composed with dilations as in (\ref{fiesa}), or they are of a very special degenerate form
that can be completely described.

As far as we know our approach is completely different from known proof techniques. It has been known before \cite{Hua}
that our problem has an equivalent reformulation as the one of characterizing coherency preserving maps on $2$ by $2$ hermitian matrices.
Combining elementary linear algebra and geometric techniques with topological tools (we will use one basic result from the homotopy theory of spheres)
we will show that continuous coherency preservers on $2$ by $2$ hermitian matrices are either degenerate, or they preserve adjacency.
It has been known before (here, nontrivial methods from geometry are needed) that non-degenerate adjacency preservers are
of the standard form. First, let us make clear what we mean by a degenerate coherency preserver:

\begin{definition}
A continuous coherency preserver $\phi : \Rea^4 \rightarrow \Rea^4$ is called \textbf{degenerate} when
its range is included in a light cone.
\end{definition}

Now, we can state our main result:

\begin{theorem}\label{glsvni}
Let $\phi : \Rea^4 \to \Rea^4$ be a continuous map satisfying (\ref{clepet}). Then either $\phi$ is a Poincar\' e similarity, or $\phi$ is degenerate.
\end{theorem}

A constructive description of all degenerate continuous coherency preservers is given in Section \ref{degeneratesection}.

\subsection{Restatement in terms of hermitian matrices}

Our problem can be put into the broader context of real quadratic spaces:
let $(V,q)$ be a real quadratic space, i.e.\ $V$ is a real vector space and $q$
is a quadratic form on $V$. The polar form $(x,y) \mapsto \frac{1}{2}(q(x+y)-q(x)-q(y))$ of $q$
is denoted by $b_q$. The quadratic space $(V,q)$ is called a \textbf{Minkowski space} when $V$ is finite-dimensional
and $q$ has signature $(1,n-1)$, where $n:=\dim V$ is greater than $1$ (in that case, $q$ is non-degenerate).

In that context, two points $x,y$ in $V$ are called \textbf{coherent} whenever $q(x-y)=0$, and they are called \textbf{adjacent} whenever
they are coherent and $x \neq y$. Given $a \in V$, the \textbf{coherency cone} of $a$ is defined as
$$\calC(a):=\{m \in V : \; q(m-a)=0\}.$$
Note that $\calC(a)=a+\calC(0)$.

A \textbf{coherency preserver} is a map $\phi : V \rightarrow V$ that satisfies
$$\forall (a,b)\in V^2, \; q(b-a)=0 \Rightarrow q(\phi(b)-\phi(a))=0;$$
such a map is called \textbf{degenerate} when the range of $\phi$ is included in a coherency cone.

A \textbf{linear similarity} of $(V,q)$ is a bijective linear transformation $u$ for which there is a nonzero scalar
$\lambda$ such that $\forall x \in V, \; q(u(x))=\lambda \,q(x)$.
An \textbf{affine similarity} of $(V,q)$ is the composite of a linear similarity with a translation.

It is known that the Poincar\'e similarities are the affine similarities of $(\Rea^4,r \mapsto \langle r,r\rangle)$.
Thus, our main theorem can be reformulated as follows:

\begin{theorem}\label{generalized}
Let $(V,q)$ be a $4$-dimensional Minkowski space and
$\phi : V \rightarrow V$ be a continuous coherency preserver.
Then $\phi$ is an affine similarity or it is degenerate.
\end{theorem}

As all Minkowski spaces with a given dimension are isometric, proving Theorem \ref{generalized} only requires that
we prove it for a particular $4$-dimensional Minkowski space.

An important example of such a space consists of ${\cal H}_2$, the space of all $2$ by $2$ hermitian matrices, equipped
with the quadratic form $\det$.
For any such matrix
$$M =\left[ \begin{matrix} { t+z & x+iy \cr x-iy & t-z \cr} \end{matrix} \right],
$$
we have
$$\det M=-x^2-y^2-z^2+t^2,$$
and it is then easily seen that $({\cal H}_2,\det)$ is a Minkowski space
(it is isometric to the standard Minkowski space $(\Rea^4,r \mapsto \langle r,r\rangle)$
through $(x,y,z,t) \mapsto \left[ \begin{matrix} { t+z & x+iy \cr x-iy & t-z \cr} \end{matrix} \right]$).
Two matrices $A$ and $B$ of ${\cal H}_2$ are coherent if and only if $A-B$ is singular,
and they are adjacent if and only if $A-B$ has rank $1$. Note that this corresponds to the usual notion of adjacency between
matrices with the same format. A map $\phi : {\cal H}_2 \to
{\cal H}_2$ is called an adjacency preserver if $\phi (A)$ and $\phi (B)$ are adjacent whenever $A$ and $B$ are adjacent. Since $A$ and $B$ are
coherent if and only if $A=B$ or $A$ and $B$ are adjacent, the problem of describing the general form of coherency preservers is closely related
to the one of characterizing adjacency preservers. The study of adjacency preservers was initiated by Hua seven decades ago. The optimal version of
his theorem for hermitian matrices has been obtained only very recently \cite{HuS}. It reads as follows.

\begin{theorem}\label{adjac}
Let $\phi : {\cal H}_2 \to {\cal H}_2$ be an adjacency preserver.
Then one of the following holds.
\begin{enumerate}
\item There exist $S \in {\cal H}_2$, a rank one orthogonal projection $R$, and
a function
$\rho : {\cal H}_2 \to \Rea$ such that $\forall A \in {\cal H}_2, \; \phi (A) = \rho (A) R + S$.
\item There exist $c \in \{ -1, 1 \}$,
an invertible $2$ by $2$ complex matrix $T$, and $S\in {\cal H}_2$ such that either
\begin{equation}\label{jj1}
\forall A \in {\cal H}_2, \; \phi (A) = c\,T A\, T^* + S,
\end{equation}
or
\begin{equation}\label{jj2}
\forall A \in {\cal H}_2, \; \phi (A) = c\,TA^t\, T^* + S.
\end{equation}
\end{enumerate}
\end{theorem}

Clearly, the map $A\mapsto \rho (A)R + S$ preserves adjacency if and only
if we have $\rho (A) \not= \rho (B)$ whenever $A$ and $B$ are adjacent. An example of such a function $\rho$ is
$\rho (A) = {\rm trace}\, (A)$. Indeed, if $A$ and $B$ are adjacent, then $B-A$ is of rank one and as it is diagonalizable
we find ${\rm trace}\, (B-A) \neq 0$.

Every map of the form (\ref{jj1}) or of the form (\ref{jj2}) will be called a \textbf{standard adjacency preserver}.

In fact, in \cite{HuS} adjacency preservers from the set of all $m$ by $m$ hermitian matrices into
the set of all $n$ by $n$ hermitian matrices, where $m$ and $n$ are arbitrary positive integers, were described.
We have formulated the result only for the 2 by 2 matrices as only this low-dimensional case is needed
to solve our problem.

Clearly, every adjacency preserver is a coherency preserver.
The first impression might be that the problems of characterizing adjacency preservers and coherency preservers
are almost the same. But it turns out that the second one is much more difficult. A first evidence might be that degenerate adjacency preserving maps are of an extremely simple form, while the structure of non-standard coherency preservers is more complicated. We have a complete description
of adjacency preservers between matrix spaces of arbitrary dimensions while at present the description of coherency
preservers seems to be out of reach even in the simplest $2$ by $2$ case.

It is well-known that the standard adjacency preservers as described in (\ref{jj1}) and (\ref{jj2}) are the affine similarities
of $({\cal H}_2, \det)$. It follows that in order to establish Theorem \ref{generalized} it suffices to prove the following result:

\begin{theorem}\label{air380}
Let $(V,q)$ be a Minkowski space whose dimension is greater than $3$, and let
$\phi : V \to V$ be a continuous coherency preserver. Then either $\phi$ preserves adjacency, or
it is degenerate.
\end{theorem}

\subsection{Structure of the article}

In Section \ref{toolssection}, we state a few basic results on Minkowski spaces that will be used throughout the rest of the article.
Section \ref{degeneratesection} features a constructive view of continuous degenerate coherency preservers on $4$-dimensional Minkowski spaces.
The last section, which is independent on the results of Section \ref{degeneratesection}, consists of the proof of Theorem \ref{air380}.
This proof is mostly self-contained, with the notable exception of the use of some tools from the homotopy theory of spheres.

\section{Basic lemmas}\label{toolssection}

\subsection{Maximal coherent sets}

The following results are well-known but we reprove them for the sake of completeness.

\begin{lemma}\label{coherentline}
Let $(V,q)$ be a Minkowski space. Let $a,b,c$ be pairwise coherent points of $V$.
Then, $a,b,c$ belong to the same affine line.
\end{lemma}

\noindent
{\sl Proof.}
Setting $x:=b-a$ and $y:=c-b$, we have $q(x)=q(y)=q(x+y)=0$. Thus, $b_q(x,y)=0$.
If $x$ and $y$ were linearly independent, $\Span(x,y)$ would be a two-dimensional
totally isotropic subspace of $V$, contradicting the fact that the Witt index of $q$ equals $1$.
\enp

A \textbf{coherent set} in a Minkowski space $(V,q)$ is a subset of $V$ whose elements are pairwise coherent.
As an obvious corollary to the above lemma, we get:

\begin{corollary}
The maximal coherent subsets of a Minkowski space $V$ are the lines of the form $a+\Rea x$ in which $a\in V$ and $x \in \calC(0) \setminus \{0\}$.
\end{corollary}

With a similar way of reasoning as above, one shows that every line that is included in $\calC(a)$ must go through $a$.

\subsection{Intersection of coherency cones}

Throughout the rest of the section, $(V,q)$ denotes an arbitrary Minkowski space.

\begin{lemma}\label{closecoherent}
Let $\|-\|$ be an arbitrary norm on the Minkowski space $V$.
There exists a constant $M>0$ such that, for all $x$ and $y$ in $V$, there exists $z \in \calC(x) \cap \calC(y)$ such that
$$\|z-x\| \leq M\,\|y-x\|.$$
\end{lemma}

\noindent
{\sl Proof.}
Denote by $n$ the dimension of $V$.
No generality is lost in assuming that $V=\Rea^n$, $q : (x_1,\dots,x_n) \mapsto x_n^2-\sum_{k=1}^{n-1} x_k^2$,
and the norm under consideration is the standard Euclidean norm $(x_1,\dots,x_n) \mapsto \sqrt{\sum_{k=1}^n x_k^2}$.
Letting $(x,y)\in V^2$, we prove that there is some $z \in \calC(x) \cap \calC(y)$ such that $\|z-x\| \leq \|y-x\|$.
In this respect, one sees that by applying a translation followed with a transformation of the form $(\overrightarrow{u},t) \mapsto
(Q\overrightarrow{u},t)$ with $Q \in O(\Rea^{n-1})$, no generality is lost in assuming that $x=(0,\dots,0)$ and $y=(0,\dots,0,a,b)$
for some $(a,b)\in \Rea^2$, in which case one checks that $z:=\left(0,\dots,0,\frac{a+b}{2},\frac{a+b}{2}\right)$ has the required properties.
\enp

%

The following result is probably known but we have not been able to find any reference to it in the literature.

\begin{lemma}\label{tri}
Let $(V,q)$ be a Minkowski space whose dimension is greater than $3$.
Let $a,b,c$ be pairwise non-collinear vectors in $\calC(0)$. Then there exists $d \in V \setminus \calC(0)$
that is coherent with $a$, $b$ and $c$.
\end{lemma}

\noindent

{\sl Proof.}
Denote by $B$ the polar form of $q$ and by $\bot$ the orthogonality relation that is associated with it.
Then, the problem amounts to finding $d \in V$ such that $q(d) \neq 0$ and $2B(a,d)=2B(b,d)=2B(c,d)=q(d)$.

The space $P:=\{a-b,b-c\}^\bot$ has dimension at least $2$.
Note that $a,b,c$ are linearly independent: indeed, the contrary would yield a $2$-dimensional subspace $W$ of $V$ that contains $a,b,c$, 
and hence the quadratic form $q$ would vanish on three distinct $1$-dimensional subspaces of $W$: it would follow that
$q$ vanishes everywhere on $W$, contradicting the fact that the Witt index of $q$ equals $1$. 
Then, $a,a-b,b-c$ are also linearly independent. 

As $q$ is non-degenerate, this leads to $P \not\subset \{a\}^\bot$;
in other words $d \mapsto B(d,a)$ does not vanish everywhere on $P$.
Moreover, as the Witt index of $q$ equals $1$ the quadratic form $q$ does not vanish everywhere on $P$ either.
Thus, the product of the polynomial maps $d \mapsto B(d,a)$ and $d \mapsto q(d)$ cannot vanish everywhere on $P$,
which yields a vector $d_0 \in P$ such that $B(d_0,a) \neq 0$ and $q(d_0)\neq 0$.

It follows that the
polynomial $q(t\,d_0-a)=t^2 q(d_0)-2 t B(d_0,a)$ in the variable $t$ has exactly two real roots: we denote by $t_0$ the nonzero one. 
Then, with $d:=t_0\,d_0$, we find that $q(d)=t_0^2\,q(d_0) \neq 0$ and $q(d-a)=0$, the latter of which leads to $q(d)=2 B(d,a)$. As $d$ belongs to $\{a-b,b-c\}^\bot$,
we have $2 B(d,a)=2 B(d,b)=2 B(d,c)$. Therefore, $d$ has the required properties.
\enp

\section{A constructive view of degenerate maps}\label{degeneratesection}

Set
$$\calQ:=\left\{(x,y,z,1) \mid (x,y,z)\in \Rea^3, \; x^2+y^2+z^2=1\right\}.$$
Note that $\calQ \subset \calC(0)$ and that every nonzero element of $\calC(0)$ can be written in a unique way as $t\,s$ with $t \in \Rea^*$ and $s \in \calQ$.
In particular, space-time events in $\calQ$ are coherent if and only if they are equal.

Using this notation, we can give a constructive description of all degenerate continuous coherency preserving maps on
$\Rea^4$:

\begin{proposition}\label{pasmozacel}
Let $\phi : \Rea^4 \to \Rea^4$ be an arbitrary map.
The following conditions are equivalent:
\begin{itemize}
\item[(a)] $\phi$ is a degenerate continuous coherency preserver.
\item[(b)] There exist a countable set $J$, a family $(U_j)_{j \in J}$ of pairwise disjoint open subsets of $\Rea^4$,
a family $(s_j)_{j \in J}$ of vectors of $\calQ$, a space-time event $s'$ and a continuous map $f : \Rea^4 \rightarrow \Rea$ such that:
\begin{itemize}
\item[(i)] For all distinct $j$ and $j'$ in $J$ and all $r \in U_j$ and $r' \in U_{j'}$, one has $\langle r-r',r-r'\rangle \neq 0$;
\item[(ii)] The map $f$ vanishes everywhere on $\Rea^4 \setminus \bigcup_{j \in J} U_j$;
\item[(iii)] For all $r \in \Rea^4$,
$$\phi(r)=\left\{\begin{matrix}
{s'+f(r).s_j  & \textrm{if}\; r \in U_j \cr
s' & \textrm{otherwise.}}
\end{matrix}\right.$$
\end{itemize}
\end{itemize}
\end{proposition}

In (b), note that in the case when $J$ is the empty set the map $\phi$ is a constant map sending all space-time events to $s'$.

\bigskip
\noindent
{\sl Proof.} Assume first that condition (b) holds. Obviously, the range of $\phi$ is included in $\calC(s')$ and hence $\phi$
is degenerate. Next, we show that $\phi$ is continuous and that it satisfies (\ref{clepet}). If the sequence of space-time events
$(r_n)$ converges to $r$ and $r\in U_j$ for some $j$, then $r_n \in U_j$ for all $n$ that are large enough, and hence $\phi (r_n)$ converges to $\phi(r)$ because $f$ is continuous.
If, on the other hand, $r$ does not belong to any of the sets $U_j$, then $f(r_n) \to f(r) = 0$ as $n$ tends to infinity, and because
$\calQ$ is bounded we see that $\phi (r_n) \to s' = \phi (r)$ as $n \to \infty$. Therefore, $\phi$ is continuous.
Let $r_1$ and $r_2$ in $\Rea^4$ be such that $ \langle \phi (r_1 ) - \phi ( r_2 ) , \phi ( r_1 ) - \phi ( r_2 ) \rangle \not=0$.
As $\forall j \in J, \; \langle s_j , s_j \rangle = 0$, we deduce that $r_1 \in U_j$ and $r_2 \in U_{j'}$ for some distinct $j$ and $j'$ in $J$
and it follows from our assumptions that $\langle r_1 -r_2 , r_1 - r_2 \rangle \not= 0$. Therefore, $\phi$ preserves coherency.

Conversely, assume that condition (a) holds. The range of $\phi$ is then included in $\calC(s')$ for some space-time event $s'$.
We can write $\phi(r)-s'=(?,?,?,f(r))$ for all $r \in \Rea^4$. The map $f : \Rea^4 \to \Rea$
is then continuous and hence $U:=f^{-1} (\Rea^*)$ is an open subset of $\Rea^4$ (possibly empty).
The map
$$g : r \in U \mapsto \frac{1}{f(r)}\,(\phi(r)-s')$$
is obviously continuous and its range is included in $\calQ$.
For $s \in \calQ$, set $U_s:=g^{-1} (\{s\})$, so that the subsets $U_s$ are pairwise disjoint and their union equals $U$.

Next, we prove that each subspace $U_s$ is open in $U$. To do so, we prove that $g$ is locally constant on $U$.
Let $r_0 \in U$. By Lemma \ref{closecoherent}, we can find a neighborhood $U'$ of $r_0$ in $U$
such that for all $r \in U'$, there exists $r' \in U$ that is coherent with both $r_0$ and $r$.
Let $r \in U'$.  Then, we obtain $r' \in U \cap \calC(r_0) \cap \calC(r)$. The space-time events
$\phi(r_0)$ and $\phi(r')$ are coherent and belong to $\calC(s') \setminus \{s'\}$.
By Lemma \ref{coherentline}, it follows that $g(r_0)$ and $g(r')$ are collinear, and hence $g(r_0)=g(r')$.
Similarly, we obtain $g(r)=g(r')$, and hence $g(r_0)=g(r)$. Thus, $g$ is locally constant on $U$.
It follows that $U_s$ is an open subset of $U$, and hence of $\Rea^4$, for all $s \in \calQ$.

Finally, set $J:=\{s \in \calQ : \; U_s \neq \emptyset\}$. As $\Rea^4$ is separable the set $J$ is countable.
Then, we know that $(U_s)_{s \in J}$ is an open cover of $U$ and that $f$ vanishes everywhere on $\Rea^4 \setminus U$,
and it is straightforward to check that
$$\forall r \in \Rea^4, \; \phi(r)=
\left\{\begin{matrix}
{s'+f(r).s  & \textrm{if}\; r \in U_s \cr
s' & \textrm{otherwise.}}
\end{matrix}\right.$$
Finally, let $j$ and $j'$ be different elements of $J$, and let
$r \in U_j$ and $r' \in U_{j'}$. Then, as $g(r)\neq g(r')$, the events $\phi(r)$ and $\phi(r')$ are non-coherent, and hence $r$ and $r'$ are non-coherent.
Thus, all the requirements in (b) are fulfilled.
\enp

 We next show that it is easy to construct such
 degenerate maps. All we need to do is to find a family of open sets $U_n$, $n \in \Nat$, and a continuous function $f$ with the above properties.
Let $U_n$, $n \in \Nat$, be an $\varepsilon$-neighbourhood (with respect to the usual Euclidean distance) of the space-time event $(n,0,0,0)$.
By choosing $\varepsilon$ small enough we find that $\langle r -s , r - s \rangle \not= 0$  for every pair of space-time events $(r,s)\in U_m \times U_n$
with distinct non-negative integers $m$ and $n$. 
If we choose arbitrary nonzero continuous functions $f_n : \Rea^4 \to \Rea$ with support in $U_n$, for $n \in \Nat$,
then it is straightforward to check that the function
$$f : r \mapsto
\left\{\begin{matrix}
{ f_n(r) & \textrm{if $r \in U_n$} \cr
0 & \textrm{otherwise}}
\end{matrix}\right.$$
is continuous.

\section{From coherency preservers to adjacency preservers}

\subsection{General considerations}

In proving Theorem \ref{air380} for all $n$-dimensional Minkowski spaces, it suffices to address
the situation of a particular one. Thus, we fix an integer $n \geq 4$ and we consider the space $V=\Rea^n$
equipped with the quadratic form
$$q : (x_1,\dots,x_n) \mapsto x_n^2-\sum_{k=1}^{n-1} x_k^2.$$
The polar form of $q$ is denoted by $B$.
In $V$, we have a particular subset
$$\calQ=\left\{(x_1,\dots,x_{n-1},1) \mid (x_1,\dots,x_{n-1}) \in \Rea^{n-1}, \; \sum_{k=1}^{n-1} x_k^2=1\right\}$$
and the linear form
$$\eta : (x_1,\dots,x_n) \mapsto x_n.$$
Every \emph{nonzero} vector $x$ of $\calC(0)$ then splits in a unique fashion as $x=t\,p$ for some $t \in \Rea^*$ and some $p \in \calQ$,
and more precisely we have $t=\eta(x)$ and $p=\frac{1}{\eta(x)}\,x$.
For any $a \in V$, we define
$$\pi_a : \left\{\begin{matrix}{ \calC(a) \setminus \{a\} & \rightarrow & \calQ \cr
m & \mapsto &  \frac{1}{\eta(m-a)}\,(m-a)}
\end{matrix}\right.$$
and we note that it is an open mapping.
Obviously, $\calQ$ is homeomorphic to the $(n-2)$-dimensional sphere $S^{n-2}$. In particular, it is compact and connected.

Clearly, if $\phi : V \to V$ is a coherency preserver, then $\phi$ maps every coherent line into some coherent line. Using
$\pi_a$, the set of all coherent lines through the point $a$ can be identified with $\calQ$.

We finish with the simple observation -- that will be used repeatedly in the rest of the section -- that
composing a map $\phi : V \to V$ with translations on both sides affects neither the assumptions nor the conclusion of Theorem \ref{air380}.

\subsection{Images of coherent lines}

\begin{proposition}\label{pikci}
Let $\phi : V \to V$ be a continuous coherency preserver.
Assume that $\phi$ is constant on some coherent line. Then $\phi$ is degenerate.
\end{proposition}

\noindent
{\sl Proof.}
No generality is lost in assuming that $\phi$ is constant on $\Rea a$ for some $a \in \calC(0) \setminus \{0\}$.
Let $x \in V$ be such that $B(a,x) \neq 0$. As $\forall t \in \Rea, \; q(ta-x)=-2 B(a,x) t+q(x)$,
some point $b \in \Rea a$ is coherent with $x$, so that $\phi(x)$ is coherent with $\phi(b)=\phi(0)$.

The set $\{x \in V : \; B(a,x)=0\}$ is a linear hyperplane of $V$ and hence
$\left\{x \in V : \; B(a,x) \neq 0\right\}$ is dense in $V$.
Since $\calC(\phi(0))$ is a closed subset of $V$ and $\phi$ is continuous, we conclude that the range of $\phi$
is included in $\calC(\phi(0))$.
\enp

We will assume from now on that the continuous coherency preserver $\phi : V \rightarrow V$ is not degenerate.
Then, by Proposition \ref{pikci}, for every point $a \in V$ and every vector
$p \in \calQ$, there exists a unique $p' \in \calQ$ such that $\phi$ maps $a+\Rea p$ into $\phi(a)+\Rea p'$.
Hence, for each $a \in V$ the map $\phi$ induces a map
$$
\varphi_a : {\cal Q} \to {\cal Q}
$$
such that
$$\forall a \in V, \; \forall p \in \calQ, \; \phi(a+\Rea p) \subset \phi(a)+\Rea \varphi_a(p).$$

\begin{lemma}\label{pinoc}
The map $\Phi : V \times {\cal Q} \to {\cal Q}$ defined by
$$
\Phi : (a,p) \mapsto  \varphi_a (p)
$$
is continuous.
\end{lemma}

\noindent
{\sl Proof.} Assume that a sequence $(a_n)\in V^{\Nat}$ converges to $a$ and  a sequence $(p_n) \in \calQ^{\Nat}$ converges to $p$. We choose
a real number $t_0$ such that $\phi (a + t_0\, p) = \phi (a) + \lambda\, \varphi_a (p)$ for some nonzero real number $\lambda$.
By the continuity of $\phi$ we know that
$$
\phi (a_n + t_0\, p_n)=\phi (a_n) + s_n\, \varphi_{a_n} (p_n)  \to \phi (a) + \lambda\, \varphi_a (p)
$$
as $n \to \infty$, and since $\phi (a_n)$ converges to $\phi (a)$, we have
$$
s_n\, \varphi_{a_n} (p_n) \to \lambda\, \varphi_a (p)
$$
as $n\to \infty$. By taking the last coordinate, it follows that $\lim s_n = \lambda$, and as $\lambda \neq 0$ we conclude that
$\lim \varphi_{a_n} (p_n) = \varphi_a (p)$, as desired.
\enp

 As $\calQ$ is compact, the compact-open topology on $C({\cal Q}, {\cal Q})$ -- the space of all continuous functions from ${\cal Q}$ to itself -- coincides with the topology of uniform convergence. A straightforward consequence of the above lemma is that the map from $V$ to $C({\cal Q}, {\cal Q})$ given by $a\mapsto \varphi_a$ is continuous.

\subsection{The existence of a non-degenerate generic point}

Let $\phi : V \to V$ be a non-degenerate continuous coherency preserver.
We denote by ${\cal W} \subset V$ the set of all points $b \in V$ such that $\phi^{-1} (\{ b \})$ has nonempty interior in $V$.
By separability of $V$, the set ${\cal W}$ is countable. We will call a point $a\in V$ \textbf{generic} if $\phi (a) \not\in {\cal W}$.

\begin{proposition}\label{existnondeg}
Let $\phi : V \to V$ be a non-degenerate continuous coherency preserver.
Then there exists a generic point $a$ such that $\varphi_a$ is nonconstant.
\end{proposition}

\noindent
{\sl Proof.}
We perform a \emph{reductio ad absurdum} by assuming that $\varphi_a$ is constant for every generic point $a$.
There are three main steps: we shall successively prove:
\begin{itemize}
\item That $\varphi_a$ is constant for all $a \in V$;
\item That $a \mapsto \varphi_a$ is constant on $V$;
\item And finally that $\phi$ maps $V$ into a coherent line.
\end{itemize}

Let $a \in V$. We shall prove that $\varphi_a$ is locally constant.
Let $p \in {\cal Q}$. As $\phi$ is continuous and nonconstant on the line $a+\Rea p$, it takes uncountably many values on it,
and hence we can find $b \in a+\Rea p$ that is generic and such that $\phi(b) \neq \phi(a)$.

Let us consider an arbitrary norm $\|-\|$ on $V$, and denote by $M \geq 1$ a constant satisfying the conclusion of Lemma \ref{closecoherent} for that norm.
As $\phi$ is continuous at $b$, we can choose $\varepsilon>0$ such that $\phi(m) \neq \phi(a)$ for all $m \in V$ such that
$\|m-b\| \leq \varepsilon$.
The set $\calC(a) \setminus \{a\}$
is homeomorphic to $\Rea^* \times \calQ$, and both $\Rea^*$ and $\calQ$ are locally connected.
It follows that there exists a connected
neighborhood $\Delta$ of $b$ in $\calC(a) \setminus \{a\}$ that lies within the closed ball with center $b$ and radius $\frac{\varepsilon}{M}\cdot$
Then, $\widetilde{\Delta}:=\pi_a(\Delta)$ is a connected neighborhood of $p$ in ${\cal Q}$ (remember that $\pi_a$ is an open mapping):
in the next step we shall prove that $\varphi_a$ is constant on $\widetilde{\Delta}$.

Let $b'$ be a generic point in $\Delta$. As $\|b'-b\| \leq \frac{\varepsilon}{M}$, we find some $m \in V$
that is coherent with both $b$ and $b'$ and that satisfies $\|m-b\| \leq M \|b'-b\| \leq \varepsilon$.
Thus, $\phi(m)$ is different from $\phi(a)$, and it must belong to both lines $(\phi(a)\phi(b))$ and $(\phi(a)\phi(b'))$ as $\varphi_b$ and $\varphi_{b'}$
are constant. It follows that $(\phi(a)\phi(b))=(\phi(a)\phi(b'))$.
We conclude that $\varphi_a(\pi_a(b'))=\varphi_a(\pi_a(b))$ for every generic point $b'$ in $\Delta$.
However, the set of all vectors $\varphi_a(\pi_a(b'))=\pi_{\phi(a)}(\phi(b'))$, with \emph{non-generic} $b' \in \Delta$, is countable.
It follows that $\varphi_a(\widetilde{\Delta})$ is countable.
On the other hand, as $\varphi_a$ is continuous and $\widetilde{\Delta}$ is connected the space $\varphi_a(\widetilde{\Delta})$
is connected. We conclude that $\varphi_a(\widetilde{\Delta})$ consists of a single point, which proves our claim.

Thus, $\varphi_a$ is locally constant. As ${\cal Q}$ is connected and $\varphi_a$ is continuous we conclude that $\varphi_a$ is constant.

Now, we move on to our second step. Let $a$ and $b$ be two points of $V$.
Assume first that $a$ and $b$ are adjacent. Denoting by
$p$ the sole vector in ${\cal Q}$ such that $\phi((ab)) \subset \phi(a)+\Rea p$, we see that both maps $\varphi_a$ and $\varphi_b$ take the value $p$,
and as they are constant we deduce that they are equal.
In the general case, Lemma \ref{closecoherent} yields a point $c$ that is coherent with both $a$ and $b$, and we deduce that $\varphi_a=\varphi_c=\varphi_b$.

Now, we have a vector $p \in {\cal Q}$ such that $\phi(b)\in \phi(a)+\Rea p$ for every coherent pair $(a,b)$.
Let $a \in V$. Then there is a point $b$ that is coherent with both $0$ and $a$, and hence there are scalars
$\alpha$ and $\beta$ such that $\phi(a)-\phi(b)=\alpha\,p$ and $\phi(b)-\phi(0)=\beta\,p$.
Thus, $\phi(a)=\phi(0)+(\alpha+\beta)\,p$. Therefore, the range of $\phi$ is included in the coherent line $\phi(0)+\Rea p$,
which is itself included in $\calC(\phi(0))$. This contradicts the assumption that $\phi$ is non-degenerate.
\enp

\subsection{When $0$ is a non-degenerate generic point}

The aim of this section is to prove the following statement, which is the main part of our proof of Theorem \ref{air380}.

\begin{proposition}\label{cinaj}
Let $\phi : V \to V$ be a continuous coherency preserver. Assume that $\phi$ is not degenerate, $\phi (0) = 0$ and $\varphi_0$ is nonconstant.
Assume further that $0$ is generic. Then $\phi$ preserves adjacency.
\end{proposition}

The proof will be carried out in a few steps that will be formulated as lemmas.

\begin{lemma}\label{cinaj1}
Let $\phi : V \to V$ be as in Proposition \ref{cinaj}. Then
$\phi$ maps $\calC(0) \setminus \{0\}$ into itself.
\end{lemma}

\noindent
{\sl Proof.} Assume on the contrary that $\phi (a) = 0$ for some $a \in \calC(0) \setminus \{0\}$.
Since $\varphi_0$ is continuous and nonconstant, its range
contains at least three different points. It follows that there exist vectors $b,c$ in $\calC(0) \setminus \{0\}$ such that $a,b,c$ are pairwise non-coherent
and $\phi (b)$ and $\phi(c)$ are non-coherent elements of $\calC(0)$. According to Lemma \ref{tri} there
exists $d \in V \setminus \calC(0)$ that is coherent with $a$, $b$, and $c$.
We claim that $\phi (m) = 0$ for all $m$ in some neighborhood of $d$, contradicting the assumption that $0$ is generic.

Among the points of $\Rea b$, we see that $b$ is coherent with $d$ but $0$ is not.
Thus, the affine map $t \mapsto q(d-tb)=q(d)-2 B(d,b) t$ is nonconstant, and hence $B(d,b)\neq 0$.
Similarly $B(d,c) \neq 0$.
Thus, as $B$ is continuous we can find a neighborhood $U$ of $d$ such that
$$\forall m \in U, \; B(m,b)B(m,c) \neq 0.$$
Let $m \in U$. One checks that the sole point of $\Rea b$ that is coherent with $m$ is
$\frac{q(m)}{2 B(m,b)}\,b$, and ditto for $b$ replaced with $c$.
As $\phi$ is continuous and $\phi(b) \neq 0$ and $\phi(c) \neq 0$, we deduce that
there is a neighborhood $U'$ of $d$ such that $U' \subset U$ and every $m \in U'$
is coherent with some $b' \in \Rea b$ for which $\phi(b') \neq 0$ and with some $c' \in \Rea c$ for which $\phi(c') \neq 0$.

Let $m \in U'$ be such that $\phi(m)$ is coherent with $0$, and let $b'$ and $c'$ be as above with respect to $m$.
Then, $\phi(m)$, $\phi(0)=0$ and $\phi(b')$ are pairwise coherent, and hence $\phi(m) \in \Rea \phi(b')$.
Similarly $\phi(m) \in \Rea \phi(c')$ and hence $\phi(m)=0$ as $\phi(b')$ and $\phi(c')$ are non-collinear.
In particular $\phi(d)=0$ since $\phi(d)$ is coherent with $\phi(a)=0$.

Finally, Lemma \ref{closecoherent} yields a neighborhood $U''$ of $d$ such that $U'' \subset U'$ and
for every $m \in U''$, there exists $m' \in U'$ that is coherent with both $d$ and $m$.
Let then $m \in U''$, and take $m'$ as above. As $\phi(m')$ is coherent with $\phi(d)=0$,
we find that $\phi(m')=0$, and then applying the same principle once more leads to $\phi(m)=0$.
Thus, $\phi$ vanishes everywhere on $U''$, contradicting the assumption that $0$ is generic.
\enp

\begin{lemma}\label{cinaj2}
Let $\phi : V\to V$ be as in Proposition \ref{cinaj}. Then the map $\varphi_0$ is injective.
\end{lemma}

\noindent
{\sl Proof.} Assume on the contrary that there exist non-collinear
vectors $a$ and $b$ in $\calC(0) \setminus \{0\}$ such that $\phi (a)$ and $\phi (b)$ are coherent.
The map $t \mapsto \phi (t\,b)$ is continuous from the unit interval and takes different values $\phi(0)$ and $\phi(b)$
on its boundary points. Hence, replacing $b$ by $t\,b$ for a suitable nonzero $t$ from the unit interval,
we may assume further that $\phi (a) \not= \phi (b)$.

The assumption that $\varphi_0$ is nonconstant yields the existence of some $c \in \calC(0) \setminus \{0\}$
such that $\phi (c)$ is not coherent with $\phi (a)$.
Then clearly, $c$ is coherent neither with $a$, nor with $b$. Lemma \ref{tri} yields some $d \in V \setminus \calC(0)$ that is
coherent with $a$, $b$, and $c$. But then $\phi (d)$ is coherent with $\phi (a)$ and $\phi (b)$, and thus $\phi (d)$ belongs to the coherent line
through these two points, which goes through $0$. The coherency of $\phi (d)$ with $\phi (c)$ then yields that
$\phi (d) = 0$.

As in the previous lemma we can find a neighborhood $U$ of $d$ such that for every $m \in U$ the point $\phi (m)$ is coherent with
$t \,\phi (a)$ and $s\,\phi (a)$
for some nonzero real numbers $t,s$ with $t\not=s$, and $\phi(m)$ is coherent also to $r\,\phi (c)$ for some nonzero  real number $r$.
But then, as before, $\phi (m) =0$ for all $m \in U$, contradicting our assumption that $0$ is generic.
\enp

\begin{lemma}\label{cinaj3}
Let $\phi : V \to V$ be as in Proposition \ref{cinaj}. Then the map $\varphi_0$ is a homeomorphism.
\end{lemma}

\noindent
{\sl Proof.} The map $\varphi_0$ is injective and continuous  from ${\cal Q}$ to itself. By the invariance of domain theorem
$\varphi_0 ({\cal Q})$ is open in ${\cal Q}$. On the other hand $\varphi_0 ({\cal Q})$ is compact in ${\cal Q}$.
As $\calQ$ is connected we deduce that $\varphi_0 ({\cal Q})={\cal Q}$.
\enp

\begin{lemma}\label{cinaj4}
Let $\phi : V \to V$ be as in Proposition \ref{cinaj}. Then for every $a\in V$ the map $\varphi_a$ is nonconstant. 
\end{lemma}

\noindent
{\sl Proof.} Here we will apply the homotopy theory of spheres.
The range of the continuous map $a \in V \mapsto \varphi_a$ is arcwise connected in $C({\cal Q}, {\cal Q})$.
Since $\varphi_0$ is a homeomorphism, every map in its arcwise component in $C({\cal Q}, {\cal Q})$ is of the same degree which is either $1$ or
$-1$, and hence such a map must be nonconstant. 
\enp

\begin{lemma}\label{cinaj5}
Let $\phi : V \to V$ be as in Proposition \ref{cinaj}. Then every $a\in V$ is generic for $\phi$.
\end{lemma}

\noindent
{\sl Proof.} Assume this is not true. Then there exists a nonempty open subset $U \subset V$ on which $\phi$ is constant. Let ${\cal D}$
be an arbitrary coherent line through some point $a$ of $U$.
We know that $\phi$ is continuous and nonconstant on ${\cal D}$ and hence the range of the restriction of
$\phi$ to  ${\cal D}$ is uncountable. Therefore we can find $b \in {\cal D}$ such that $\phi(b) \not\in {\cal W}$,  the set of all points  in $V$ whose $\phi$-preimages have nonempty interior in $V$.
In particular, $a \not=b$ and $b$ is generic. We have $a-b = s\,p$ for some nonzero real number $s$ and some vector $p \in \calQ$.
If $p'\not=p$ is a vector of $\calQ$ close enough to $p$,
then $b + s\,p' \in U$. Both lines, one going through $b$ and $a$, and the other one going through $b$ and $b+s\,p'$ are mapped into the line going through $\phi (b)$ and $\phi (a) = \phi (b+s\,p')$.
Thus, the map $\varphi_b$ is not injective. This contradicts Lemma \ref{cinaj2} applied to the map $m \mapsto \phi (m+b) -\phi (b)$.
\enp

We are now ready to complete the proof of the main result of this section.

\bigskip

\noindent
{\sl Proof of Proposition \ref{cinaj}.}
Let $a$ and $b$ be adjacent points of $V$. Then $b$ is generic and $\varphi_b$ is nonconstant. Applying Lemma \ref{cinaj1} to the map  $m \mapsto \phi (m+b) -\phi (b)$, we conclude that
$\phi ( (a-b) + b) - \phi (b)$ is a nonzero element of $\calC(0)$, as desired.
\enp

\subsection{Conclusion}

Now, we can finish the proof of Theorem \ref{air380}. Let
$\phi : V \to V$ be a continuous coherency preserver.
Assume that $\phi$ is non-degenerate. By Proposition \ref{existnondeg}, there exists a generic point $a$ such that $\varphi_a$ is nonconstant.
Then, the continuous coherency preserver $m \mapsto \phi(m+a)-\phi(m)$ satisfies the assumptions of Proposition \ref{cinaj}, and
hence it preserves adjacency. Therefore, $\phi$ preserves adjacency.

\end{document}